\documentclass[11pt,a4paper]{article}

\usepackage{theorem,enumerate}
\usepackage{amsmath,latexsym,amssymb,amsfonts}
\usepackage{fancybox}
\usepackage{graphicx}
\usepackage[all]{xy}

\newcounter{Scounter}
\theorembodyfont{\normalfont\slshape}
\newtheorem{thm}{Theorem}

\newtheorem{Thm}{Theorem}

\newtheorem{prop}[thm]{Proposition}
\newtheorem{lem}{Lemma}

\newtheorem{claim}{Claim}[section]

\newtheorem{fact}[claim]{Fact} 

\newtheorem{Con}[Thm]{Conjecture}

\newtheorem{pro}[thm]{Problem}
\theoremstyle{remark}
\newtheorem{rema}{Remark}
{\theorembodyfont{\upshape}
\newtheorem{rem}[rema]{Remark}}

\newcommand{\proof}{\medbreak\noindent\textit{Proof.}\quad}

\newcommand{\qed}{{$\quad\square$\vs{3.6}}}

\newcommand{\vs}[1]{\vspace*{#1 mm}}

\numberwithin{equation}{section}

\newcommand{\ora}{\overrightarrow}
\newcommand{\ola}{\overleftarrow}

\makeatletter
\newcommand{\xRightarrow}[2][]{%
\ext@arrow 0055{\Rightarrowfill@}{#1}{#2}%
}
\def\Leftarrowfill@{\arrowfill@\Leftarrow\Relbar\Relbar}
\newcommand{\xLongleftrightarrow}[2][]{%
\ext@arrow 0055{\llrafill@}{#1}{#2}%
}
\def\llrafill@{\arrowfill@\Leftarrow\Relbar\Rightarrow}
\makeatother

\title{On directed 2-factors in digraphs and\\ 
$2$-factors containing perfect matchings in bipartite graphs}

\author{Shuya Chiba$^{1}$\thanks{This work was supported by JSPS KAKENHI grant 17K05347, 17K05348}
\thanks{E-mail address: \texttt{schiba@kumamoto-u.ac.jp}} 
\hspace{+12pt}
Tomoki Yamashita$^{2}$\thanks{This work was supported by JSPS KAKENHI grant 16K05262}
\thanks{E-mail address: \texttt{yamashita@math.kindai.ac.jp}} 
\vspace{+8pt}
 \\
\small
$^1$\small\textsl{Applied Mathematics, Faculty of Advanced Science and Technology, 
Kumamoto University,}\\ 
\vspace{+6pt}
\small\textsl{2-39-1 Kurokami, Kumamoto 860-8555, Japan}\\
\small
$^{2}$\small\textsl{Department of Mathematics, Kindai University,}\\
\small\textsl{3-4-1 Kowakae, Higashi-Osaka, Osaka 577-8502, Japan}
}

\date{}

\begin{document}
\setlength{\baselineskip}{17pt}

\maketitle

\vspace{-24pt}
\begin{abstract} 
In this paper, 
we give the following result: 
If $D$ is a digraph of order $n$,
and if $d_{D}^{+}(u) + d_{D}^{-}(v) \ge n$ 
for every two distinct vertices $u$ and $v$ with $(u, v) \notin A(D)$, 
then $D$ has a directed $2$-factor with exactly $k$ directed cycles of length at least $3$, 
where $n \ge 12k+3$. 
This result is equivalent to the following result: 
If $G$ is a balanced bipartite graph of order $2n$ with partite sets $X$ and $Y$, 
and if $d_{G}(x)+d_{G}(y) \ge n + 2$
for every two vertices $x \in X$ and $y \in Y$ with $xy \notin E(G)$, 
then 
for every perfect matching $M$, 
$G$ has a $2$-factor with exactly $k$ cycles of length at least $6$ 
containing every edge of $M$, 
where $n \ge 12k+3$. 
These results are generalizations of 
theorems concerning Hamilton cycles due to Woodall (1972) and Las Vergnas (1972), 
respectively.


\end{abstract}

\section{Introduction}
\label{sec:introduction}

We consider only finite graphs. 
For standard graph-theoretic terminology 
not explained in this paper,
we refer the reader to \cite{DiestelBook}. 
Unless stated otherwise, ``graph'' means a simple undirected graph. 
Let $G$ be a graph. 
We denote by $V(G)$ and $E(G)$ 
the vertex set and the edge set of $G$, respectively. 
For a vertex $v$ of $G$, 
$d_{G}(v)$ denotes the degree of $v$ in $G$.
Let $\sigma_2(G)$ be the minimum degree sum
of two non-adjacent vertices in $G$, i.e., 
$\sigma_2(G)=\min \big\{ d_G(u) + d_G(v): u, v \in V(G), u \neq v, uv \notin E(G) \big\}$
if $G$ is not complete; 
otherwise, $\sigma_2(G) = + \infty$.

A graph is said to be \textit{hamiltonian} if it has a \textit{Hamilton cycle}, i.e., a cycle containing all the vertices. 
The following $\sigma_{2}$ condition, due to Ore (1960), 
is classical and well known in graph theory 
and there are many results on degree conditions which generalize it 
(see a survey \cite{LiSurvey2013}).

\begin{Thm}[Ore \cite{Ore1960}]
\label{Thm:Ore1960}
Let $G$ be a graph of order $n \ge 3$. 
If $\sigma_{2}(G) \ge n$, then $G$ is hamiltonian. 
\end{Thm}

A $2$-\textit{factor} of a graph is a spanning subgraph in which every component is a cycle, 
and thus a Hamilton cycle is a $2$-factor with exactly $1$ cycle. 
Brandt, Chen, Faudree, Gould and Lesniak (1997) gave the following $\sigma_{2}$ condition 
for the existence of a $2$-factor with exactly $k$ cycles.

\begin{Thm}[Brandt et al. \cite{BCFGL1997}]
\label{Thm:BCFGL1997}
Let $k$ be a positive integer, 
and let $G$ be a graph of order $n \ge 4k$. 
If $\sigma_{2}(G) \ge n$, then $G$ has a $2$-factor with exactly $k$ cycles. 
\end{Thm}

Theorem~\ref{Thm:Ore1960} is essentially the case $k = 1$ of Theorem~\ref{Thm:BCFGL1997}, 
since the result of Theorem~\ref{Thm:Ore1960} is easy to see when $n = 3$; 
thus Theorem~\ref{Thm:BCFGL1997} implies Theorem~\ref{Thm:Ore1960}.

For a digraph $D$, 
we denote by $A(D)$ the arc set of $D$, 
and 
let $d_{D}^{+}(v)$ and $d_{D}^{-}(v)$ be the out-degree and the in-degree of a vertex $v$ in $D$, respectively. 
A \textit{directed Hamilton cycle} is a directed cycle containing all the vertices of the digraph, 
and 
a \textit{directed $2$-factor} of a digraph is a spanning subdigraph in which every component is a directed cycle.

In \cite{Woodall1972}, 
Woodall (1972) gave the digraph version of Theorem~\ref{Thm:Ore1960} as follows.

\begin{Thm}[Woodall \cite{Woodall1972}]
\label{Thm:Woodall1972}
Let $D$ be a digraph of order $n \ge 2$.  
If $d_{D}^{+}(u) + d_{D}^{-}(v) \ge n$ 
for every two distinct vertices $u$ and $v$ with $(u, v) \notin A(D)$, 
then $D$ has a directed Hamilton cycle. 
\end{Thm}

In fact, the following remark implies that 
this theorem is a generalization of Theorem~\ref{Thm:Ore1960}.

\begin{rem}
\label{rem:D_{G}}
For a given graph $G$, 
let $D_{G}$ be the digraph of order $|V(G)|$ 
obtained from $G$ 
by replacing each edge $uv$ in $G$ with 
two arcs $(u, v)$ and $(v, u)$. 
It is easy to see that a graph $G$ satisfies the hypothesis (conclusion) of Theorem~\ref{Thm:Ore1960} 
if and only if $D_{G}$ satisfies the hypothesis (conclusion) of Theorem~\ref{Thm:Woodall1972}; 
thus Theorem~\ref{Thm:Woodall1972} implies Theorem~\ref{Thm:Ore1960}. 
\end{rem}

In this paper, 
we show that the Woodall condition also implies the existence of the following directed $2$-factor, 
which is our main result.

\begin{thm}
\label{thm: sigma1,1 for directed 2-factor with k cycles}
Let $k$ be a positive integer, 
and let $D$ be a digraph of order $n$, 
where $n \ge 12k+3$. 
If $d_{D}^{+}(u) + d_{D}^{-}(v) \ge n$ 
for every two distinct vertices $u$ and $v$ with $(u, v) \notin A(D)$, 
then $D$ has a directed $2$-factor with exactly $k$ directed cycles of length at least $3$. 
\end{thm}

\begin{figure}[h]
\begin{center}
\[
\xymatrix@C=60pt@R=3pt{
    \textup{{\small undirected graphs}}  
& \textup{{\small digraphs}} 
& \textup{{\small bipartite graphs}}
\\
    \textup{\fbox{Theorem~\ref{Thm:BCFGL1997}}} 
    \ar@{=>}[ddddd]^(.4){\textup{{\tiny ($n \ge 4$)}}} 
    \ar@{<=}[r]^{\textup{Remark~\ref{rem:D_{G}}}}_{\textup{{\tiny ($n \ge 12k + 3$)}}} 
& \textup{\doublebox{Theorem~\ref{thm: sigma1,1 for directed 2-factor with k cycles}}} 
    \ar@{=>}[ddddd]^(.4){\textup{{\tiny ($n \ge 15$)}}} 
    \ar@{<=>}[r]^{\textup{Remark~\ref{rem:from D to G}}}
& \textup{\doublebox{Theorem~\ref{thm: sigma1,1 for M-2-factor with k cycles}}} 
    \ar@{=>}[ddddd]^(.4){\textup{{\tiny ($n \ge 15$)}}}  
\\
&
&
\\
&
&
\\
&
&
\\
&
&
\\
    \textup{\fbox{Theorem~\ref{Thm:Ore1960}}} \ar@{<=}[r]^{\textup{Remark~\ref{rem:D_{G}}}} 
& \textup{\fbox{Theorem~\ref{Thm:Woodall1972}}} \ar@{<=>}[r]^{\textup{Remark~\ref{rem:from D to G}}} 
& \textup{\fbox{Theorem~\ref{Thm:Las Vergnas1972}}}  
}
\]
\end{center}
\caption{The relation between results in this paper}
\label{relation}
\end{figure}
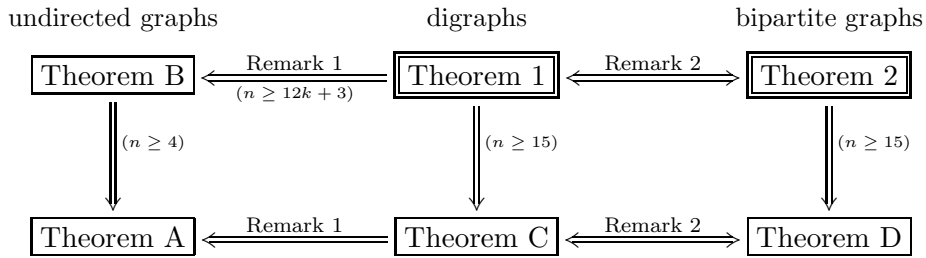

By Remark~\ref{rem:D_{G}}, 
Theorem~\ref{thm: sigma1,1 for directed 2-factor with k cycles} implies 
the result of Theorem~\ref{Thm:BCFGL1997} 
for graphs with order $n \ge 12k + 3$, 
and it also clearly implies the result of Theorem~\ref{Thm:Woodall1972} 
for graphs with order $n \ge 15$. 
Thus, in a sense, 
Theorem~\ref{thm: sigma1,1 for directed 2-factor with k cycles} 
is a common generalization of 
Theorems~\ref{Thm:BCFGL1997} and \ref{Thm:Woodall1972} (see Figure~\ref{relation}).

The degree condition in Theorem~\ref{thm: sigma1,1 for directed 2-factor with k cycles} 
is best possible in the following sense. 
Consider the complete bipartite graph $G = K_{(n-1)/2, (n + 1)/2}$, 
and let $D_{G}$ be the digraph obtained from $G$ by the same procedure as 
in Remark~\ref{rem:D_{G}}. 
Then, 
$\min \big\{ d_{D_{G}}^{+}(u) + d_{D_{G}}^{-}(v) : (u, v) \notin A(D_{G}), u \neq v\big\} = n-1$ 
and clearly $D_{G}$ does not have a directed $2$-factor.

On the other hand, 
the order condition in Theorem~\ref{thm: sigma1,1 for directed 2-factor with k cycles} 
comes from our proof techniques. 
The complete bipartite graph $K_{2k-1, 2k-1}$ shows that 
$n \ge 4k-1$ is necessary for the existence of a $2$-factor with 
exactly $k$ cycles in simple undirected graphs, 
and hence it follows from the similar argument as above 
that $n \ge 4k-1$ is also necessary for Theorem~\ref{thm: sigma1,1 for directed 2-factor with k cycles}. 
Unfortunately, our proof of Theorem~\ref{thm: sigma1,1 for directed 2-factor with k cycles} 
requires the stronger condition $n \ge 12k + 3$; 
this arises from the condition $n \ge 12k - 9$ 
which is needed for the final contradiction in the proof of Theorem~\ref{thm:small k disjoint M-cycles} 
(see Section~\ref{sec:proof of packing theorem}).

In the next section, 
we further give a relationship 
between 
Theorem~\ref{Thm:Woodall1972}, Theorem~\ref{thm: sigma1,1 for directed 2-factor with k cycles} and 
the results on $2$-factors containing perfect matchings in bipartite graphs
(Theorem~\ref{Thm:Las Vergnas1972} 
and Theorem~\ref{thm: sigma1,1 for M-2-factor with k cycles}) in Figure~\ref{relation}. 
The proof of Theorem~\ref{thm: sigma1,1 for directed 2-factor with k cycles} 
is presented in Sections~\ref{sec:proof of main}--\ref{sec:proof of partition theorem}.

\section{2-factors containing perfect matchings in bipartite graphs}
\label{sec:2-factors containing perfect matchings in bipartite graphs}

In previous section, 
in order to generalize Theorems~\ref{Thm:Ore1960}, \ref{Thm:BCFGL1997} and \ref{Thm:Woodall1972}, 
we have considered the directed $2$-factors with exactly $k$ directed cycles of length at least $3$  
and have given Theorem~\ref{thm: sigma1,1 for directed 2-factor with k cycles}. 
It is also known that this problems have a connection with the $2$-factor problems in bipartite graphs. 
In fact, 
Theorem~\ref{Thm:Woodall1972} is equivalent to the following theorem due to Las Vergnas (1972).  
Here, 
an edge subset $M$ of a graph $G$ is called a \textit{matching} 
if no two edges in $M$ have a common end. 
In particular, 
a matching $M$ is said to be \textit{perfect} if every vertex of $G$ is contained in some edge of $M$. 
An \textit{alternating cycle with respect to a matching $M$} 
is a cycle such that the edges belong to $M$ and not to $M$, alternatively. 
For a bipartite graph $G$ with partite sets $X$ and $Y$, 
we define 
$\sigma_{1,1}(G) = \min \big\{d_G(x)+d_G(y) : x \in X, y \in Y, xy \not\in E(G) \big\}$ 
if $G$ is not a complete bipartite graph; 
otherwise, $\sigma_{1,1}(G) = + \infty$.

\begin{Thm}[Las Vergnas \cite{Las Vergnas1972}]
\label{Thm:Las Vergnas1972}
Let $G$ be a balanced bipartite graph of order $2n \ge 4$ 
and $M$ be a perfect matching of $G$. 
If $\sigma_{1, 1}(G) \ge n + 2$, then $G$ has a Hamilton cycle containing every edge of $M$. 
\end{Thm}

\begin{figure}[h]
\begin{center}
\includegraphics[scale=1.00,clip]{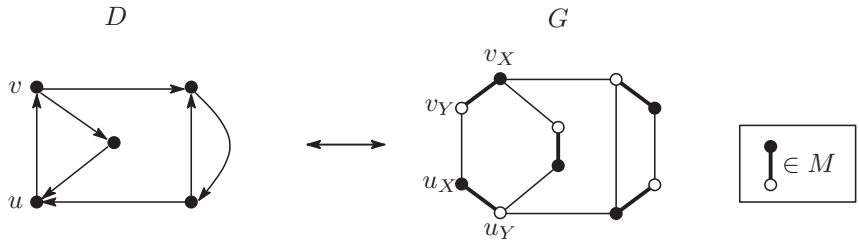}
\caption{The relationship between digraphs and bipartite graphs with a perfect matching}
\label{from bip to di}
\end{center}
\end{figure}

\begin{rem}[see also \cite{ZZW2013}]
\label{rem:from D to G}
For a given digraph $D$, 
consider the following undirected graph $G$: 
We split each vertex $v$ in $D$ into two vertices $v_{X}$ and $v_{Y}$ 
and replace each arc $(u, v)$ in $A(D)$ with a simple edge $u_{X}v_{Y}$, 
and we add the perfect matching $M = \{v_{X}v_{Y} : v \in V(D)\}$. 
Then, the resultant graph $G$ is a balanced bipartite graph of order $2|V(D)|$ 
with partite sets $\{v_{X} : v \in V(D)\}$ and $\{v_{Y} : v \in V(D)\}$ 
(see Figure~\ref{from bip to di}). 
On the other hand, 
for a given balanced bipartite graph $G$ 
with partite sets $X, Y$ and a perfect matching $M$ in $G$, 
let $D$ be the digraph of order $|V(G)|/2$ 
obtained from $G$ 
by replacing each edge $xy \in E(G) \setminus M$ ($x \in X$, $y \in Y$) 
with an arc from $x$ to $y$, 
and by contracting all edges of $M$ (see Figure~\ref{from bip to di}). 
Note that, in this construction, 
the following hold: 
\begin{itemize}
\item $(u, v) \in A(D)$ if and only if 
$u_{X}v_{Y} \in E(G)$ 
(in particular, $d_{D}^{+}(v) = d_{G}(v_{X}) - 1$ and $d_{D}^{-}(v) = d_{G}(v_{Y}) - 1$), and 
\item 
an alternating cycle of length $2l \ (\ge 4)$ with respect to $M$ in $G$ 
corresponds to a directed cycle of length $l \ ( \ge 2)$ in $D$. 
\end{itemize}
Therefore, 
there is a one-to-one correspondence between 
the class of digraphs satisfying the Woodall condition 
and 
the class of bipartite graphs satisfying the Las Vergnas condition. 
\end{rem}

This also implies that 
Theorem~\ref{thm: sigma1,1 for directed 2-factor with k cycles} is equivalent to the following theorem 
(see also Figure~\ref{relation}). 
Related results can be found 
in \cite{AFGP2007} and a survey \cite{Gould2009}.

\begin{thm}
\label{thm: sigma1,1 for M-2-factor with k cycles}
Let $k$ be a positive integer, 
and let $G$ be a balanced bipartite graph of order $2n$ 
and $M$ be a perfect matching of $G$, 
where 
$n \ge 12k + 3$. 
If $\sigma_{1,1}(G) \ge n + 2$, then $G$ has a $2$-factor with exactly $k$ cycles of length at least $6$ 
containing every edge of $M$. 
\end{thm}

\begin{prop}
\label{prop:equivalent}
Theorem~\ref{thm: sigma1,1 for directed 2-factor with k cycles} 
and Theorem~\ref{thm: sigma1,1 for M-2-factor with k cycles} are equivalent. 
\end{prop}

Therefore, 
the Las Vergnas condition also implies the existence of 
a $2$-factor 
with a prescribed number of cycles 
containing the specified perfect matching in bipartite graphs. 
In this sense, there is no difference between Hamilton cycles and 
$2$-factors with $k \ (\ge 2)$ cycles.

\section{Proof of Theorem~\ref{thm: sigma1,1 for directed 2-factor with k cycles}}
\label{sec:proof of main}

By Proposition~\ref{prop:equivalent}, 
to prove Theorem~\ref{thm: sigma1,1 for directed 2-factor with k cycles}, 
it suffices to show Theorem~\ref{thm: sigma1,1 for M-2-factor with k cycles}. 
Therefore, 
in this section, 
we introduce the steps of the proof of Theorem~\ref{thm: sigma1,1 for M-2-factor with k cycles} 
and also give two theorems in order to prove it. 
Here, 
for a bipartite graph $G$ and a matching $M$ of $G$, 
a cycle $C$ of $G$ is called an $M$-\textit{cycle} 
if $|E(C) \cap M| = \frac{|C|}{2}$, i.e., 
$C$ is an alternating cycle with respect to $M$. 
In particular, 
$C$ is called an $M$-\textit{Hamilton cycle}
if $C$ is a Hamilton cycle and an $M$-cycle,
and
a $2$-factor of a graph $G$ is called an $M$-$2$-\textit{factor}
if every component is an $M$-cycle.

The proof of Theorem~\ref{thm: sigma1,1 for M-2-factor with k cycles} 
involves two steps, 
summarized in the following two theorems.

\begin{thm}
\label{thm:small k disjoint M-cycles}
Let $k$ be a positive integer, 
and let $G$ be a balanced bipartite graph of order $2n$ 
and $M$ be a perfect matching of $G$, 
where 
$n \ge 12k - 9$. 
If $\sigma_{1,1}(G) \ge n + 2$, then $G$ contains 
$k$ disjoint $M$-cycles of length $6$ or $8$. 
\end{thm}

\begin{thm}
\label{thm:from k+1 disjoint cycles to M-2-factor with k cycles}
Let $k$ be a positive integer, 
and let $G$ be a balanced bipartite graph of order $2n > 6(k+1)$ 
and $M$ be a perfect matching of $G$. 
Suppose that $G$ contains $k+1$ disjoint $M$-cycles of length at least $6$.  
If $\sigma_{1,1}(G) \ge n + 2$, then 
 $G$ has an $M$-$2$-factor with exactly $k$ cycles of length at least $6$.
\end{thm}

Now we prove Theorem~\ref{thm: sigma1,1 for directed 2-factor with k cycles} assuming Theorems~\ref{thm:small k disjoint M-cycles}
and \ref{thm:from k+1 disjoint cycles to M-2-factor with k cycles}.

\medskip
\noindent
\textbf{Proof of Theorem~\ref{thm: sigma1,1 for directed 2-factor with k cycles}.}~We first show that 
Theorem~\ref{thm: sigma1,1 for M-2-factor with k cycles} is true. 
Let $k,~n,~G,~M$ be the same as in Theorem~\ref{thm: sigma1,1 for M-2-factor with k cycles}. 
Since 
$n \ge 12 k + 3 = 12(k+1) - 9$ and $\sigma_{1, 1}(G) \ge n + 2$, 
it follows from Theorem~\ref{thm:small k disjoint M-cycles} that 
$G$ contains $k+1$ disjoint $M$-cycles of length at least $6$. 
Then,
since $n \ge 12k + 3 > 3(k+1)$, 
it follows from Theorem~\ref{thm:from k+1 disjoint cycles to M-2-factor with k cycles}
that
$G$ has an $M$-$2$-factor consisting of $k$ cycles of length at least $6$, 
that is, Theorem~\ref{thm: sigma1,1 for M-2-factor with k cycles} is true. 
Then, by Proposition~\ref{prop:equivalent}, 
Theorem~\ref{thm: sigma1,1 for directed 2-factor with k cycles} is also true. 
\qed

Therefore, 
it suffices to show that 
Theorems~\ref{thm:small k disjoint M-cycles}
and \ref{thm:from k+1 disjoint cycles to M-2-factor with k cycles} hold. 
The proofs of Theorem~\ref{thm:small k disjoint M-cycles} 
and \ref{thm:from k+1 disjoint cycles to M-2-factor with k cycles} are presented 
in Sections~\ref{sec:proof of packing theorem} and \ref{sec:proof of partition theorem}, respectively.
In Section~\ref{sec:related problems},
we mention a problem 
related to Theorem~\ref{thm:small k disjoint M-cycles}.

\medskip
Finally, we prepare terminology and notations, 
which will be used in the proofs. 
Let $G$ be a graph. 
We denote by $N_{G}(v)$ the neighborhood of a vertex $v$ in $G$. 
For $S \subseteq V(G)$, 
let $G[S]$ denote the subgraph induced by $S$ in $G$, 
and let $G - S = G[V(G) \setminus S]$. 
For $S, T \subseteq V(G)$ with $S \cap T = \emptyset$, 
$E_{G}(S, T)$ denotes the set of edges of $G$ between $S$ and $T$, 
and let $e_{G}(S, T) = |E_{G}(S, T)|$. 
We often identify a subgraph $H$ of $G$ with its vertex set $V(H)$ 
(e.g., we often use $E_{G}(F, H)$ instead of $E_{G}(V(F), V(H))$ 
for disjoint subgraphs $F$ and $H$ of $G$). 
We denote by $P[x, y]$ a path $P$ with ends $x$ and $y$ in $G$ 
and $|P|$ denotes the number of vertices in $P$. 
Next let $G$ be a bipartite graph, 
and $M$ be a matching of $G$. 
For a subgraph $H$ of $G$, 
let $M_{H} = M \cap E(H)$. 
A path $P = P[x, y]$ 
is called an $M$-\textit{path} of $G$ 
if 
$P$ is an alternating path (i.e., a path 
such that the edges belong to $M$ and not to $M$, alternatively) 
joining $x$ and $y$ starting and ending with edges in $M$. 
In particular, 
$P$ is called an $M$-\textit{Hamilton path} 
if $P$ is also a Hamilton path of $G$. 
If $X$ and $Y$ are partite sets of $G$, 
then 
for $A \subseteq X$ (resp., $B \subseteq Y$), 
we define $\overline{A} = \{y \in Y : xy \in M$ with $x \in A\}$ 
(resp., $\overline{B} = \{x \in X : xy \in M$ with $y \in B\}$).

\section{Proof of Theorem~\ref{thm:small k disjoint M-cycles}}
\label{sec:proof of packing theorem}

In this section, we give the proof of Theorem~\ref{thm:small k disjoint M-cycles}. 
In order to 
prove it, 
we use the following two lemmas (Lemmas~\ref{lem:smaller M-cycles} and \ref{lem:degree condition for M-cycle of length 6}).

\begin{lem}
\label{lem:smaller M-cycles}
Let $G$ be a bipartite graph and $M$ be a matching of $G$, 
and let 
$C$ be an $M$-cycle and $P = P[x, y]$ be an $M$-path in $G-V(C)$. 
Assume that $e_{G}(\{x, y\}, C) > n$, where $n = \frac{|C|}{2}$. 
Then $G[V(P \cup C)]$ contains an $M$-cycle of length $|P| + 2i$ 
for each $i$, $1 \le i \le n$. 
\end{lem}
\noindent
\textbf{Proof of Lemma~\ref{lem:smaller M-cycles}.}~Let the vertices be labelled $u_{1}v_{1} \dots u_{n}v_{n}$ 
in order round $C$, 
where $u_{1}v_{1}, \dots,$ $u_{n}v_{n}$ are the edges of $M_{C}$ 
and 
$v_{1}, \dots, v_{n}$ are in the same partite set as $x$. 
If there is no $M$-cycle of length $|P| + 2i$, 
then $G$ contains at most one of the edges $xu_{j}, yv_{j + i - 1}$ 
for each $j$ ($i, j \in \{1, \dots, n\}$), 
where the subscripts are interpreted modulo $n$, 
and so $e_{G}(\{x, y\}, C) \le n$. 
This contradiction proves the result. 
\qed

\begin{lem}
\label{lem:degree condition for M-cycle of length 6}
Let $G$ be a balanced bipartite graph of 
order $2n$ 
with partite sets $X$ and $Y$, 
and 
let $M$ be a perfect matching of $G$. 
If $d_{G}(x) \ge \frac{n+3}{2}$ for every vertex $x$ in $X$, 
then 
$G$ has an $M$-cycle of length $6$. 
\end{lem}
\noindent
\textbf{Proof of Lemma~\ref{lem:degree condition for M-cycle of length 6}.}~Note that 
the hypothesis implies $n \ge \frac{n + 3}{2}$, so that $n \ge 3$. 
Also, $\sum_{y \in Y}d_{G}(y) = |E(G)| = \sum_{x \in X}d_{G}(x) \ge n (\frac{n + 3}{2})$, 
which implies that 
there exists a vertex $y$ in $Y$ 
such that 
$d_{G}(y) \ge \frac{n + 3}{2}$.

Suppose that 
$G$ contains no $M$-cycle of length $6$. 
Let $x$ be a vertex in $X$ such that $xy \in M$, 
and let $x'y' \in M \setminus \{xy\}$ 
such that $xy' \in E(G)$. 
Then $N_{G}(y) \cap ( \overline{N_{G}(x')} \setminus \{x, x'\}) = \emptyset$ 
since otherwise, 
$G$ contains an $M$-cycle of length $6$, a contradiction. 
Note that $|\overline{N_{G}(x')}| = |N_{G}(x')| = d_{G}(x') \ge \frac{n + 3}{2}$, 
ans so
\begin{align*}
d_{G}(y) \le |X| - |\overline{N_{G}(x')}| + 2 \le n - \frac{n +3}{2} + 2 = \frac{n + 1}{2} < \frac{n +3}{2}, 
\end{align*}
a contradiction. 
\qed

Now we are ready to prove Theorem~\ref{thm:small k disjoint M-cycles}.

\bigskip
\noindent
\textbf{Proof of Theorem~\ref{thm:small k disjoint M-cycles}.}~If $n = 3$, 
then by the degree condition, 
we can easily check that $G$ is a complete bipartite graph, 
and hence the assertion clearly holds. 
Thus, 
we may assume that 
$n \ge 4$, i.e., 
$n \ge \max \{4, 12 k - 9\}$.

We suppose that 
\begin{align}
\label{G has no t+1 disjoint M-cycles}
\begin{array}{cc}
\hspace{-48pt}G \textup{ contains } t \textup{ disjoint } \textup{$M$-cycles of length } 6 \textup{ or } 8 \textup{ with } 0 \le t \le k - 1, \\[1mm]
\textup{but } G \textup{ does not contain } t + 1 \textup{ disjoint } \textup{$M$-cycles of length } 6 \textup{ or } 8. 
\end{array}
\end{align}
Let $C_{1}, \dots, C_{t}$ be $t$ disjoint $M$-cycles of length $6$ or $8$, 
and 
choose $C_{1}, \dots, C_{t}$ so that $\sum_{i = 1}^{t}|C_{i}|$ is as small as possible. 
Without loss of generality, we may assume that $C_{1}, \dots, C_{t_{1}}$ are cycles of length $6$ 
and $C_{t_{1}+1}, \dots, C_{t}$ are cycles of length $8$ 
for some $t_{1}$ with $0 \le t_{1} \le t$. 
Let $H = G - \bigcup_{i=1}^{t}V(C_{i})$.

Now let 
$\mathcal{P}$ be a set of mutually disjoint $M$-paths of order $4$ in $H$, 
and we define 
\begin{align*}
\textup{$\mathcal{P}^{*} = \{P[x, y] \in \mathcal{P} : d_{G}(x) + d_{G}(y) \ge \sigma_{1, 1}(G)\}.$} 
\end{align*}
We choose $\mathcal{P}$ so that 
\begin{enumerate}[]
\item{(P1)} $|\mathcal{P}^{*}|$ is as large as possible, and 
\item{(P2)} $|\mathcal{P}|$ is as large as possible, subject to (P1). 
\end{enumerate}

\begin{claim}
\label{at least k+1 M-paths}
$|\mathcal{P}| \ge k + 1$. 
\end{claim}
\proof 
Note that $|H| \ge 2n - 8t \ge 2n - 8 (k-1) = 2n - 8k + 8$. 
Suppose that 
$|\mathcal{P}| \le k$, 
and let $H' = H - \bigcup_{P \in \mathcal{P}}V(P)$. 
Then $|H'| \ge |H| - 4k \ge 2n - 12k + 8 \ge 2 \max \{4, 12k-9\} - 12k + 8
\ge 4$, 
and so there are two distinct edges $x_{1}y_{1}$ and $x_{2}y_{2}$ in $M_{H'}$ 
(note that $M_{H'}$ is a perfect matching of $H'$), 
where $x_{1}$ and $x_{2}$ belong to the same partite set of $G$. 
Then 
$d_{H'}(x_{h}) = d_{H'}(y_{h}) = 1$ for $h \in \{1, 2\}$, 
as otherwise $H'$ has an $M$-path of order $4$, which contradicts (P2). 
In particular, $x_{1}y_{2}, x_{2}y_{1} \notin E(G)$. 
Also, 
$e_{G}(\{x_{h}, y_{h}\}, C_{i}) \le |C_{i}| = 6$ for $h \in \{1, 2\}$ and $1 \le i \le t_{1}$; 
$e_{G}(\{x_{h}, y_{h}\}, C_{i}) \le 4 < 6$ for $h \in \{1, 2\}$ and $t_{1} + 1 \le i \le t$, 
since otherwise $G[\{x_{h}, y_{h}\} \cup V(C_{i})]$ contains an $M$-cycle $C_{i}'$ of length $6$ 
by Lemma~\ref{lem:smaller M-cycles}, 
and replacing the cycle $C_{i}$ by this cycle $C_{i}'$ would violate the minimality of $\sum_{i = 1}^{t}|C_{i}|$; 
and $e_{G}(\{x_{h}, y_{h}\}, P) \le 3$ for $h \in \{1, 2\}$ and each $P \in \mathcal{P}$, 
since otherwise, it is easy to see that $G[\{x_{h}, y_{h} \} \cup V(P)]$ has an $M$-cycle of length $6$, 
which contradicts (\ref{G has no t+1 disjoint M-cycles}). 
Since $x_{1}y_{2}, x_{2}y_{1} \notin E(G)$, 
it follows that  
\begin{align*}
2n + 4 
&\le 2\sigma_{1, 1}(G) \\
&\le \big( d_{G}(x_{1}) + d_{G}(y_{2}) \big) + \big( d_{G}(x_{2}) + d_{G}(y_{1}) \big) \\
&= \big( d_{G}(x_{1}) + d_{G}(y_{1}) \big) + \big( d_{G}(x_{2}) + d_{G}(y_{2}) \big) \\
&\le 2 \big(2 +  6t + 3k \big) \\
&\le 2 \big(2 +  6(k-1) + 3k \big) = 18k - 8, 
\end{align*}
which implies $n \le 9k - 6 < \max \{4, 12k-9\}$, a contradiction. 
\qed

\begin{claim}
\label{number of M-paths with high degree}
$|\mathcal{P}^{*}| \ge |\mathcal{P}| - 1$. 
\end{claim}
\proof 
Suppose that $|\mathcal{P}^{*}| \le |\mathcal{P}| - 2$, 
and let 
$P_{1}[x_{1}, y_{1}], P_{2}[x_{2}, y_{2}] \in \mathcal{P} \setminus \mathcal{P}^{*}$ 
with $P_{1} \neq P_{2}$. 
Since $P_{1}, P_{2} \notin \mathcal{P}^{*}$, 
we have $x_{1}y_{1}, x_{2}y_{2} \in E(G)$. 
We may assume that
$x_{1}$ and $x_{2}$ belong to the same partite set of $G$. 
If $E_{G}(\{x_{1}, y_{1}\}, \{x_{2}, y_{2}\}) = \emptyset$, 
then 
\begin{align*}
\sum_{h \in \{ 1, 2\} } \big( d_{G}(x_{h}) + d_{G}(y_{h}) \big)
&= \big( d_{G}(x_{1}) + d_{G}(y_{2}) \big) + \big( d_{G}(x_{2}) + d_{G}(y_{1}) \big) \ge 2\sigma_{1, 1}(G), 
\end{align*}
and this implies that 
$d_{G}(x_{1}) + d_{G}(y_{1}) \ge \sigma_{1, 1}(G)$ 
or 
$d_{G}(x_{2}) + d_{G}(y_{2}) \ge \sigma_{1, 1}(G)$, 
which contradicts the assumption that $P_{1}, P_{2} \notin \mathcal{P}^{*}$. 
Thus $E_{G}(\{x_{1}, y_{1}\}, \{x_{2}, y_{2}\}) \neq \emptyset$. 
We will assume that $x_{1}y_{2} \in E(G)$.

Write 
$P_{1} = x_{1}y_{1}'x_{1}'y_{1}$ and $P_{2} = x_{2}y_{2}'x_{2}'y_{2}$. 
Note that 
$x_{h}y_{h}', x_{h}'y_{h} \in M$ for $h \in \{1, 2\}$. 
Consider the $M$-path $Q = y_{1}'x_{1}y_{2}x_{2}'$. 
If $y_{1}'x_{2}' \notin E(G)$, 
then $Q[y_{1}', x_{2}']$ is an $M$-path of order $4$ 
such that $d_{G}(y_{1}') + d_{G}(x_{2}') \ge \sigma_{1, 1}(G)$, 
and hence 
for $\mathcal{Q} = (\mathcal{P} \setminus \{P_{1}, P_{2}\}) \cup \{Q\}$, 
we have 
$|\mathcal{Q}^{*}| > |\mathcal{P}^{*}|$ 
because $P_{1}, P_{2} \notin \mathcal{P}^{*}$, 
which contradicts (P1). 
Thus $y_{1}'x_{2}' \in E(G)$. 
Since 
$x_{1}y_{1}, x_{1}y_{2}, y_{1}'x_{2}' \in E(G)$, 
it follows that 
$y_{1}'x_{2}, x_{1}'y_{2} \notin E(G)$; 
otherwise, 
$x_{2}y_{2}'x_{2}'y_{2}x_{1}y_{1}'x_{2}$ 
or 
$x_{1}'y_{1}x_{1}y_{1}'x_{2}'y_{2}x_{1}'$ 
is an $M$-cycle of length $6$, which contradicts (\ref{G has no t+1 disjoint M-cycles}). 
Therefore, 
\begin{align*}
&\hspace{+12pt}\big( d_{G}(x_{1}') + d_{G}(y_{1}') \big) + \big( d_{G}(x_{2}) + d_{G}(y_{2}) \big) \\
&= \big( d_{G}(x_{1}') + d_{G}(y_{2}) \big) + \big( d_{G}(y_{1}') + d_{G}(x_{2}) \big) \ge 2\sigma_{1, 1}(G), 
\end{align*}
and this implies that 
$d_{G}(x_{1}') + d_{G}(y_{1}') \ge \sigma_{1, 1}(G)$ 
or 
$d_{G}(x_{2}) + d_{G}(y_{2}) \ge \sigma_{1, 1}(G)$. 
Since $P_{2} \in \mathcal{P} \setminus \mathcal{P}^{*}$, 
we have $d_{G}(x_{1}') + d_{G}(y_{1}') \ge \sigma_{1, 1}(G)$.

Now consider the $M$-path $R_{1} = x_{1}'y_{1}x_{1}y_{1}'$, 
and let $\mathcal{R} = (\mathcal{P} \setminus \{P_{1}\}) \cup \{R_{1}\}$. 
Since 
$P_{1}[x_{1}, y_{1}] \in \mathcal{P} \setminus \mathcal{P}^{*}$ 
and 
$R_{1}[x_{1}', y_{1}'] \in \mathcal{R}^{*}$, 
we have $|\mathcal{R}^{*}| > |\mathcal{P}^{*}|$, 
which contradicts (P1) again. 
\qed

\begin{claim}
\label{|V(C_{i})|/2 + 1}
There exist at least two distinct paths 
$P[x, y]$ in $\mathcal{P}$ 
such that 
$e_{G}(\{x, y\}, C_{i}) \le 4$ 
for $1 \le i \le t_{1}$.
\end{claim}
\proof 
Suppose not. 
Then, 
for every path $P[x, y]$ in $\mathcal{P}$, except at most one, 
there exists a cycle $C_{i}$ with $1 \le i \le t_{1}$ 
such that $e_{G}(\{x, y\}, C_{i}) \ge 5$. 
Since $|\mathcal{P}| \ge k + 1$ by Claim~\ref{at least k+1 M-paths}, 
and since $t_{1} \le k - 1$, 
it follows from the Pigeonhole Principle that 
there exist two distinct paths $P[x, y]$, $P'[x', y']$ in $\mathcal{P}$ 
and a cycle $C_{i}$ with $1 \le i \le t_{1}$
such that 
$e_{G}(\{x, y\}, C_{i}) \ge 5$ and $e_{G}(\{x', y'\}, C_{i}) \ge 5$. 
Hence we can take two distinct edges 
$uv, u'v'$ in $M_{C_{i}}$ 
such that $e_{G}(\{x, y\}, \{u, v\}) = e_{G}( \{x', y'\}, \{u', v'\}) = 2$, 
and then 
it is easy to check that 
$G[V(P \cup P') \cup \{u, v, u', v'\}]$ contains two disjoint $M$-cycles of length $6$, 
which contradicts (\ref{G has no t+1 disjoint M-cycles}). 
\qed

By Claims~\ref{number of M-paths with high degree} and \ref{|V(C_{i})|/2 + 1}, 
there exists an $M$-path $P_{0}[x_{0}, y_{0}]$ in $\mathcal{P}$ 
such that 
$P_{0} \in \mathcal{P}^{*}$ 
and 
$e_{G}(\{x_{0}, y_{0}\}, C_{i}) \le 4$ 
for $1 \le i \le t_{1}$. 
Note that $e_{G}(\{x_{0}, y_{0}\}, C_{i}) \le \frac{|C_{i}|}{2} = 4$ for $t_{1} + 1 \le i \le t$ 
as well, as otherwise 
$G[V(P_{0} \cup C_{i})]$ contains an $M$-cycle of length $|P_{0}| + 2 = 6$ by Lemma~\ref{lem:smaller M-cycles}, 
and replacing $C_{i}$ by this new cycle would violate the minimality of $\sum_{i=1}^{t}|C_{i}|$. 
Let $H' = H - V(P_{0})$. Then 
\begin{align}
\label{number of edges x_{0}, y_{0} and H'}
e_{G}(\{x_{0}, y_{0}\}, H') 
&
\ge 
\sigma_{1, 1}(G) - \sum_{i = 1}^{t}e_{G}(\{x_{0}, y_{0}\}, C_{i}) - d_{G[V(P_{0})]}(x_{0}) - d_{G[V(P_{0})]}(y_{0}) \notag \\
&
\ge 
(n + 2) - 4t -2 - 2 \notag \\
& = 
n - 4t - 2. 
\end{align}

We let 
\begin{align}
\label{definition of X, Y}
\begin{array}{ll}
Y_{x_{0}} = N_{H}(x_{0}) \cap V(H'), &\quad X_{x_{0}} = \overline{Y_{x_{0}}},  \\[1mm]
X_{y_{0}} = N_{H}(y_{0}) \cap V(H'), &\quad Y_{y_{0}} = \overline{X_{y_{0}}}. 
\end{array}
\end{align}

\begin{figure}[h]
\begin{center}
\includegraphics[scale=1.00,clip]{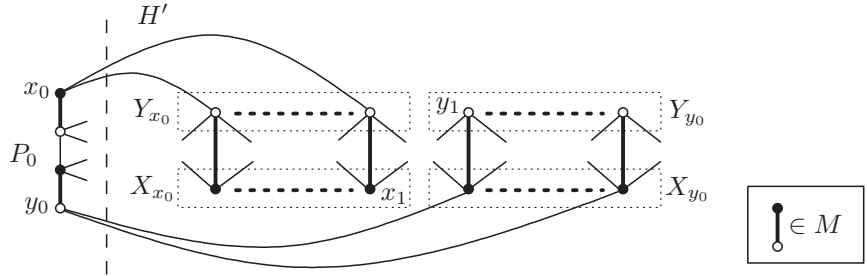}
\caption{Vertex subsets $X_{x_{0}}, Y_{x_{0}}, X_{y_{0}}$ and $Y_{y_{0}}$}
\label{subsetsX-Y}
\end{center}
\end{figure}

Since 
$H$ does not contain an $M$-cycle of length $6$ by (\ref{G has no t+1 disjoint M-cycles}), 
we have $X_{x_{0}} \cap X_{y_{0}} = Y_{x_{0}} \cap Y_{y_{0}} =  \emptyset$, 
and hence 
\begin{align}
\label{pairwise disjoint}
\textup{$X_{x_{0}}, Y_{x_{0}}, X_{y_{0}}$ and $Y_{y_{0}}$ are pairwise disjoint.}
\end{align}
Note also that 
\begin{align}
\label{same cardinality}
e_{G}(x_{0}, H') = |X_{x_{0}}| = |Y_{x_{0}}| 
\textup{ and }
e_{G}(y_{0}, H') = |X_{y_{0}}| = |Y_{y_{0}}|. 
\end{align}

\noindent
Let 
$G_{x_{0}} = G[X_{x_{0}} \cup Y_{x_{0}}]$
and
$G_{y_{0}} = G[X_{y_{0}} \cup Y_{y_{0}}]$.
Then 
by (\ref{definition of X, Y})--(\ref{same cardinality}), 
it follows that 
$G_{x_{0}}$ and $G_{y_{0}}$ are disjoint balanced bipartite graphs
with a perfect matching whose edges belong to $M$. 
We further define 
\begin{align}
\label{definition of n_{i}}
\begin{array}{ll}
&\displaystyle n_{1} = \sum_{i=1}^{t}\frac{|C_{i}|}{2}, \quad
\displaystyle n_{2} = \frac{|G_{x_{0}}|}{2} = e_{G}(x_{0}, H'), \quad 
\displaystyle n_{3} = \frac{|G_{y_{0}}|}{2} = e_{G}(y_{0}, H'), \\[1mm]
&\displaystyle n_{4} =  2 + \frac{|H'|}{2} -  \big( \frac{|G_{x_{0}}|}{2}+\frac{|G_{y_{0}}|}{2} \big). 
\end{array}
\end{align}

\begin{claim}
\label{there exists a x_{1} with deg <= n_{2}/2+1}
If $n_{2} > 0$ $($resp., $n_{3} > 0)$, then 
there exists a vertex $x_{1}$ in $X_{x_{0}}$ 
$($resp., a vertex $y_{1}$ in $Y_{y_{0}})$ 
such that 
$d_{G_{x_{0}}}(x_{1}) \le \frac{n_{2} + 2}{2}$
$\big($resp., 
$d_{G_{y_{0}}}(y_{1}) \le \frac{n_{3} + 2}{2}\big)$. 
\end{claim}
\proof 
Suppose that $n_{2} > 0$ 
and that 
$d_{G_{x_{0}}}(x) \ge \frac{n_{2} + 3}{2}$ 
for every vertex $x$ in $X_{x_{0}}$. 
Since
$M_{G_{x_{0}}} = M \cap E(G_{x_{0}})$ is a perfect matching of $G_{x_{0}}$ 
and 
$d_{G_{x_{0}}}(x) \ge \frac{n_{2} + 3}{2}$ 
for $x \in X_{x_{0}}$, 
it follows from Lemma~\ref{lem:degree condition for M-cycle of length 6}
that
$G_{x_{0}}$ has an $M$-cycle of length $6$, 
which contradicts (\ref{G has no t+1 disjoint M-cycles}).

By symmetry, if $n_{3} > 0$, then 
there exists a vertex $y_{1}$ in $Y_{y_{0}}$ 
such that 
$d_{G_{y_{0}}}(y_{1}) \le \frac{n_{3} + 2}{2}$. 
\qed

By (\ref{number of edges x_{0}, y_{0} and H'}) and (\ref{definition of n_{i}}), 
and since $t \le k - 1$ and $n \ge 12k-9$, 
we have 
\begin{align}
\label{lower bound of n_{2} + n_{3}}
n_{2} + n_{3} 
= 
e_{G}(\{x_{0}, y_{0}\}, H')  \ge  n - 4t - 2 \ge 8k - 7 > 0. 
\end{align}
In particular, 
without loss of generality, 
we may assume that 
$n_{2} > 0$. 
Then by Claim~\ref{there exists a x_{1} with deg <= n_{2}/2+1}, 
there exists a vertex 
$x_{1}$ in $X_{x_{0}}$ 
such that 
$d_{G_{x_{0}}}(x_{1}) \le \frac{n_{2} + 2}{2}$. 
If $n_{3} > 0$, then let $y_{1}$ be also the vertex as in Claim~\ref{there exists a x_{1} with deg <= n_{2}/2+1}; 
otherwise, let $y_{1} = y_{0}$. 
Note that 
\begin{align}
\label{x_{1}y_{1} is not an edge}
N_{G}(x_{1}) \cap (\{y_{0}\}) \cup Y_{y_{0}}) = \emptyset. 
\end{align}
since otherwise, 
$H$ has an $M$-cycle of length $6$ or $8$, which contradicts (\ref{G has no t+1 disjoint M-cycles}).

\begin{claim}
\label{degree of x_{1} and y_{1}}
We have 
$d_{G}(x_{1}) \le n_{1} + \frac{n_{2}}{2} + n_{4}$ 
and 
$d_{G}(y_{1}) \le n_{1} + \frac{n_{3}}{2} + n_{4}$. 
\end{claim}
\proof 
We first show that 
$d_{G}(x_{1}) \le n_{1} + \frac{n_{2}}{2} + n_{4}$. 
Note that $e_{G}(x_{1}, P_{0}) \le 1$ since $x_{1}y_{1} \notin E(G)$ by (\ref{x_{1}y_{1} is not an edge}). 
Combining this with (\ref{definition of X, Y})--(\ref{definition of n_{i}}), (\ref{x_{1}y_{1} is not an edge}) 
and the definition of $x_{1}$, we get 
\begin{align*}
d_{G}(x_{1}) 
&
= 
\sum_{i=1}^{t}e_{G}(x_{1}, C_{i}) 
+ e_{G}(x_{1}, P_{0})
+ d_{H'}(x_{1}) \\
&
\le 
n_{1} + 1 
+ \Big( 
d_{G_{x_{0}}}(x_{1}) 
+ 
|N_{H'}(x_{1}) \cap Y_{y_{0}}| 
+ \big( \frac{|H'|}{2} - |Y_{x_{0}}| - |Y_{y_{0}}| \big) \Big) \\
&
\le 
n_{1} + 1 
+ \Big( 
\frac{n_{2} + 2}{2}
+ 
0 
+ \big( \frac{|H'|}{2} - \frac{|G_{x_{0}}|}{2}-\frac{|G_{y_{0}}|}{2} \big) \Big) \\
&
= 
n_{1} + 
\frac{n_{2}}{2} 
+ 
\Big(2 + \frac{|H'|}{2} - \frac{|G_{x_{0}}|}{2}-\frac{|G_{y_{0}}|}{2} \Big) 
= 
n_{1} + 
\frac{n_{2}}{2} + n_{4}. 
\end{align*}

We next show that 
$d_{G}(y_{1}) \le n_{1} + \frac{n_{3}}{2} + n_{4}$. 
If $n_{3} > 0$, 
then this holds by the same argument. 
Thus we may assume $n_{3} = 0$. 
Recall that, in this case, $y_{1} = y_{0}$, 
and note that 
by (\ref{definition of n_{i}}), 
$e_{G}(y_{1}, H') = n_{3} = 0$ 
and 
$n_{4} \ge 2$. 
Hence, 
\begin{align*}
d_{G}(y_{1}) 
= 
\sum_{i=1}^{t}e_{G}(y_{1}, C_{i}) 
+ d_{G[V(P_{0})]}(y_{1}) 
\le 
n_{1} + 2 
\le
n_{1} + \frac{n_{3}}{2} + n_{4}. 
\quad\quad\square
\end{align*}
\vs{3.6}

Since $x_{1}y_{1} \notin E(G)$ by (\ref{x_{1}y_{1} is not an edge}), 
it follows from Claim~\ref{degree of x_{1} and y_{1}} 
and the hypothesis of Theorem~\ref{thm:small k disjoint M-cycles} that 
\begin{align*}
n + 2 
\le \sigma_{1, 1}(G) 
\le d_{G}(x_{1}) + d_{G}(y_{1}) 
\le 2n_{1} + \frac{n_{2} + n_{3}}{2} + 2n_{4}. 
\end{align*}
This inequality implies that $2(n + 2) \le 4n_{1} + n_{2} + n_{3} + 4n_{4}$. 
Since $n=n_{1} + n_{2} + n_{3} + n_{4}$ by (\ref{definition of n_{i}}) 
and the fact that $|H| = |H'| + |P_{0}| = |H'| + 4$, 
it follows that 
$n \le 3(n_{1} + n_{4}) - 4$. 
By (\ref{lower bound of n_{2} + n_{3}}), 
we obtain
$n_{1} + n_{4} = n - (n_{2} + n_{3}) \le 4t+2$. 
Thus
we have 
\begin{align*}
n \le 3(4t+2) - 4 = 12t +2 \le 12k - 10 < 12k - 9, 
\end{align*}
a contradiction.

This completes the proof of Theorem~\ref{thm:small k disjoint M-cycles}. 
\qed

\section{Proof of Theorem~\ref{thm:from k+1 disjoint cycles to M-2-factor with k cycles}}
\label{sec:proof of partition theorem}

We first prepare terminology and notation.
Let $G$ be a graph. 
We write a cycle (or a path) $C$ with a given orientation as $\ora{C}$. 
Let $\ora{C}$ be an oriented cycle (or path). 
For $x \in V(C)$, 
we denote the successor and the predecessor of $x$ on $\ora{C}$ 
by $x^{+}$ and $x^{-}$. 
For $X \subseteq V(C)$, 
let 
$X^{-} = \{x^{-} : x \in V(C)\}$. 
For $x, y \in V(C)$,
we denote by $x\ora{C}y$
the path with ends $x$ and $y$ on $\ora{C}$.
The reverse sequence of $x\ora{C}y$ is denoted by $y\ola{C}x$. 
In the rest of this paper, 
we 
assume that 
every cycle has a fixed orientation.

\medskip

In order to prove Theorem~\ref{thm:from k+1 disjoint cycles to M-2-factor with k cycles}, 
we give three lemmas as follows (Lemmas~\ref{lem:insertion edge}--\ref{lem:M-2-factor}). 
(We omit the proof of Lemma~\ref{lem:insertion edge} 
since 
it is easy 
and 
it is similar to the proof of Lemma~\ref{lem:smaller M-cycles} 
in Section~\ref{sec:proof of packing theorem}).

\begin{lem}
\label{lem:insertion edge}
Let $G$ be a bipartite graph and $M$ be a matching of $G$, 
$C$ be an $M$-cycle or an $M$-path, 
and $x$ and $y$ be two vertices of $G-V(C)$ 
which belong to different partite sets of $G$. 
\begin{enumerate}[{\upshape(1)}]
\item
\label{insertion edge in cycle}
If $C$ is a cycle and $e_{G}(\{x, y\}, C) \ge \frac{|C|}{2} + 1$, 
then there exists an edge $uv$ in $E(C) \setminus M$ such that 
$e_{G}(\{x, y\}, \{u, v\}) = 2$. 

\item
\label{insertion edge in path}
If $C$ is a path and $e_{G}(\{x, y\}, C) \ge \frac{|C|}{2} + 2$, 
then there exists an edge $uv$ in $E(C) \setminus M$ such that 
$e_{G}(\{x, y\}, \{u, v\}) = 2$. 
\end{enumerate}
\end{lem}

Before stating Lemmas~\ref{lem:crossing} and \ref{lem:M-2-factor}, 
we further prepare the following terminology. 
Let $G$ be a bipartite graph and $M$ be a matching of $G$. 
For an $M$-cycle $C$ (resp., an $M$-path $C$) and an $M$-path $Q[x, y]$ 
such that $V(C) \cap V(Q) = \emptyset$, 
$Q$ is \textit{insertible} in $C$ if 
there exists an edge $uv$ in $E(C) \setminus M$ such that 
$e_{G}(\{x, y\}, \{u, v\}) = 2$. 
We call the edge $uv$ an \textit{insertion edge} of $Q$. 
Note that by Lemma~\ref{lem:insertion edge}, 
if $C$ is a cycle and $e_{G}(\{x, y\}, C) \ge \frac{|C|}{2} + 1$, 
then 
$Q$ is insertible in $C$; 
if $C$ is a path and $e_{G}(\{x, y\}, C) \ge \frac{|C|}{2} + 2$, 
then 
$Q$ is insertible in $C$. 
Note also that if $C$ is a cycle and $Q$ is insertible in $C$,  
then 
$G[V(C \cup Q)]$ has an $M_{C \cup Q}$-Hamilton cycle; 
if $C$ is a path and $Q$ is insertible in $C$, 
then 
$G[V(C \cup Q)]$ has an $M_{C \cup Q}$-Hamilton path 
such that the ends of the path are the same ends as $C$. 
An $M$-path $P$ is \textit{maximal} if there exists no $M$-path $Q$ such that $|P|<|Q|$ and $V(P)  \subseteq V(Q)$.

\begin{lem}
\label{lem:crossing}
Let $G$ be a balanced bipartite graph of order $2n$, 
$M$ be a perfect matching of $G$,
and
$D_{1}, \dots, D_{s}$ be $s$ 
disjoint $M$-cycles 
of length at least $6$ in $G$. 
Let $H = G- \bigcup_{i=1}^{s}V(D_{i})$,
and
$\ora{P_{0}}[x, y]$ be 
a maximal $M$-path of order at least $2$ in $H$. 
If $\sigma_{1, 1}(G) \ge n + 2$, 
then one of the following (\ref{P0})--(\ref{degree sum of x and y in outside}) holds: 
\begin{enumerate}[{\upshape(i)}]

\item 
\label{P0}
either $|P_{0}| = 2$ or $G[V(P_{0})]$ has an $M_{P_{0}}$-Hamilton cycle, 

\item 
\label{P0 cup Di}
for some cycle $D_{i}$ with $1 \le i \le s$, 
$G[V(P_{0} \cup D_{i})]$ has an $M_{P_{0} \cup D_{i}}$-2-factor 
with exactly two cycles $D_{0}$ and $D_{i}'$ of length at least $6$ such that $V(D_{i}) \subsetneq V(D_{i}')$, 

\item
\label{degree sum of x and y in outside}
$e_{G}(\{x,y\}, G-H) \ge |G-H|/2 + 1$.

\end{enumerate}
\end{lem}
\noindent
\textbf{Proof of Lemma~\ref{lem:crossing}.}~Suppose that 
each of (\ref{P0}) and (\ref{P0 cup Di}) does not hold; we show that (\ref{degree sum of x and y in outside}) holds. 
Since (\ref{P0}) does not hold, 
we have 
$xy \notin E(G)$, 
and hence 
$d_{G}(x) + d_{G}(y) \ge n + 2$. 
Since ${P_{0}}$ is maximal,
we obtain $e_{G}(\{x,y\}, H-P_{0})=0$. 
Therefore, 
it suffices to show that $d_{G[V(P_{0})]}(x) + d_{G[V(P_{0})]}(y) \le \frac{|P_{0}|}{2} + 1$, 
since this will imply that 
\begin{align*}
e_{G}(\{x, y\}, G - H) \ge (n + 2) - \big( \frac{|P_{0}|}{2} + 1 \big) \ge \frac{|G-H|}{2} + 1. 
\end{align*}

\begin{figure}[h]
\begin{center}
\includegraphics[scale=1.00,clip]{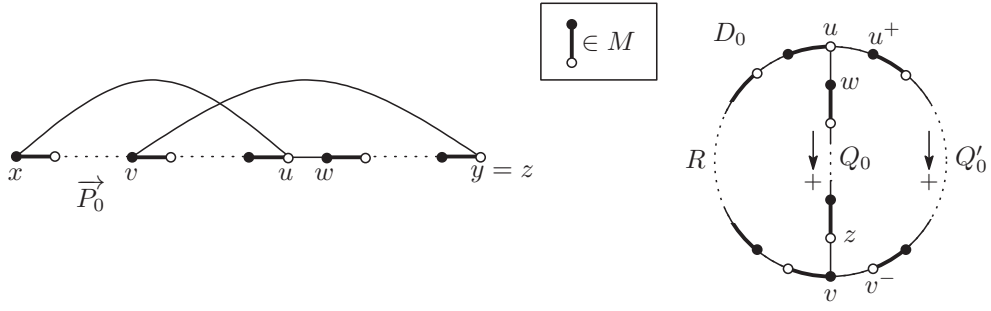}
\caption{The $M$-cycle $D_{0}$ and the $M$-path $Q_{0}$}
\label{pathPprime}
\end{center}
\end{figure}

Suppose that $d_{G[V(P_{0})]}(x) + d_{G[V(P_{0})]}(y) \ge \frac{|P_{0}|}{2} + 2$. 
Let $v$ be the first vertex along $\ora{P_{0}}$ that is in $N_{G}(y)$, 
and let $u$ be the last vertex along $\ora{P_{0}}$ that is in $N_{G}(x)$. 
Since all vertices of $N_{G[V(P_{0})]}(x)^{-}$ and $N_{G[V(P_{0})]}(y)$ are in the same partite set, 
it follows that there are at most $|P_{0}|/2$ vertices in the union of these two sets 
and hence at least two vertices in their intersection. 
It follows from this and the fact that 
$xy \notin E(G)$ that 
the vertices $x, v, u, y$ are distinct and occur in this order along $\ora{P_{0}}$, and $u \neq v^{+}$. 
Let $\ora{D_{0}}$ be the $M$-cycle $x \ora{P_{0}} u x$ and let $\ora{Q_{0}}[w, z]$ be the $M$-path $u^{+} \ora{P_{0}}y$ 
(see the left figure of Figure~\ref{pathPprime}). 
Then the following hold: 
\begin{enumerate}[]
\item{\rm (A)}
$V(D_{0}) \cap V(Q_{0}) = \emptyset$ and $V(D_{0}) \cup V(Q_{0}) = V(P_{0})$; 

\item{\rm (B)}
there exist two vertices $u, v \in V(D_{0})$ such that $uv \notin M$ and $wu, zv \in E(G)$; 

\item{\rm (C)}
$|D_{0}| \ge 6$. 
\end{enumerate}

In the rest of this proof, 
we use only (A)--(C). 
We make no further use of the orientation of $\ora{P_{0}}$; 
the superscripts $^{+}$ and $^{-}$ will refer to the orientations of $\ora{D_{0}}$ and $\ora{Q_{0}}$. 
We assume $\ora{D_{0}}$ is oriented so that $uu^{+} \in E(D_{0}) \setminus M$. 
Note that $u^{+} \neq v$ since (\ref{P0}) does not hold. 
Choose $D_{0}, Q_{0}$ and $u, v$ so that (A)--(C) hold and 
\begin{enumerate}[]
\item{\rm(D1)} $|D_{0}|$ is as large as possible, and 
\item{\rm(D2)} $|u\ora{D_{0}}v|$ is as small as possible, subject to (D1). 
\end{enumerate}

\noindent
Let $D_{0}' = v \ora{D_{0}} u w \ora{Q_{0}} zv$ and $Q_{0}' = u^{+} \ora{D_{0}} v^{-}$, 
and let $R = v \ora{D_{0}} u$ (note that $D_{0}' = Q_{0} \cup R$). 
See the right figure of Figure~\ref{pathPprime}.

Suppose that $N_{G}(w) \cap V(Q_{0}') \neq \emptyset$, 
say $b \in N_{G}(w) \cap V(Q_{0}')$. 
Then $|b \ora{D_{0}} v| < |u \ora{D_{0}} v|$ and $bv \notin M$. 
Replacing the pair of vertices $(u, v)$ with $(b, v)$, 
this contradicts (D2). 
Thus $N_{G}(w) \cap V(Q_{0}') = \emptyset$. 
Similarly, we have $N_{G}(z) \cap V(Q_{0}') = \emptyset$.

Suppose next that $N_{G}(u^{+}) \cap V(Q_{0}) \neq \emptyset$, 
say $b \in N_{G}(u^{+}) \cap V(Q_{0})$. 
Consider the $M$-cycle $D = w \ora{Q_{0}} b u^{+} \ora{D_{0}} u w$. 
If $b = z$, then 
(\ref{P0}) holds, contradicting our assumption. 
Thus $b \neq z$, 
and we further consider the $M$-path $Q = b^{+} \ora{Q_{0}} z$. 
Then $V(D) \cap V(Q) = \emptyset$ 
and $V(D) \cup V(Q) = V(P_{0})$. 
Furthermore, 
$|D| > |D_{0}| \ge 6$ 
and 
$b^{+}b, zv$ are two independent edges with $b, v \in V(D)$ and $bv \notin M$. 
This contradicts (D1). 
Thus $N_{G}(u^{+}) \cap V(Q_{0}) = \emptyset$. 
By the similar way, 
we see that $N_{G}(v^{-}) \cap V(Q_{0}) = \emptyset$.

By the above arguments, 
we can get 
\begin{align}
\label{deg sum of w and z in Q0', u+ and v- in Q0}
e_{G}(\{w, z\}, Q_{0}') = e_{G}(\{u^{+}, v^{-}\}, Q_{0}) = 0. 
\end{align}
Moreover, 
we clearly have 
\begin{align}
\label{deg sum of w and z in Q0, u+ and v- in Q0'}
d_{G[V(Q_{0})]}(w) + d_{G[V(Q_{0})]}(z) \le |Q_{0}| 
\textup{ and }
d_{G[V(Q_{0}')]}(u^{+}) + d_{G[V(Q_{0}')]}(v^{-}) \le |Q_{0}'|. 
\end{align}

On the other hand, 
since (\ref{P0}) does not hold, 
it follows that 
$Q_{0}$ is not insertible in $D_{0}$. 
This in particular implies that 
$Q_{0}$ is not insertible in the $M$-path $R$. 
Similarly, $Q_{0}'$ is not insertible in $D_{0}'$, 
and hence $Q_{0}'$ is not also insertible in the $M$-path $R$. 
Therefore, by Lemma~\ref{lem:insertion edge}(\ref{insertion edge in path}), we get 
\begin{align}
\label{deg sum of u+ and v- in R, w and z in R}
e_{G}(\{w, z\}, R) \le \frac{|R|}{2}+1 \textup{ and } e_{G}(\{u^{+}, v^{-}\}, R) \le \frac{|R|}{2}+1. 
\end{align}

Recall that by (C), $|D_{0}| \ge 6$. 
Since (\ref{P0 cup Di}) does not hold, 
it follows that $Q_{0}$ is not insertible in each $D_{i}$ ($1 \le i \le s$). 
Moreover, since $uv \notin M$ by (B), we have $|R| \ge 4$, that is, $|D_{0}'| = |Q_{0}| + |R| \ge 6$. 
This implies that 
$Q_{0}'$ is not also insertible in each $D_{i}$ ($1 \le i \le s$). 
Therefore, by Lemma~\ref{lem:insertion edge}(\ref{insertion edge in cycle}), we get 
\begin{align}
\label{deg sum of w ,z u+ and v- in D_{i}}
e_{G}(\{w, z, u^{+}, v^{-}\}, D_{i}) \le |D_{i}| \textup{ for } 1 \le i \le s.  
\end{align}
Since ${P_{0}}$ is maximal,
it follows that
\begin{align}
\label{deg sum of w ,z u+ and v- in H-P_{0}}
e_{G}(\{w, z, u^{+}, v^{-}\}, H-P_{0})=0.
\end{align}

By (\ref{deg sum of w and z in Q0', u+ and v- in Q0})--(\ref{deg sum of w ,z u+ and v- in H-P_{0}}), 
we have 
\begin{align*}
2n + 4 
&
\le d_{G}(w) + d_{G}(z) + d_{G}(u^{+}) + d_{G}(v^{-}) \\
& 
\le |Q_{0}| + |Q_{0}'| + 2 \times \Big( \frac{|R|}{2} + 1 \Big) + \sum_{i=1}^{s}|D_{i}| 
\le |G| + 2 
= 2n + 2,
\end{align*}
a contradiction. 

This contradiction shows that $d_{G[V(P_{0})]}(x) + d_{G[V(P_{0})]}(y) \le \frac{|P_{0}|}{2} + 1$, 
which, as explained in the first paragraph, completes the proof of Lemma~\ref{lem:crossing}. 
\qed

\begin{lem}
\label{lem:M-2-factor}
Let $k,~n,~G,~M$ be the same as in Theorem~\ref{thm:from k+1 disjoint cycles to M-2-factor with k cycles}. 
Under the same degree sum condition as Theorem~\ref{thm:from k+1 disjoint cycles to M-2-factor with k cycles}, 
$G$ has an $M$-$2$-factor with 
exactly $k+1$ or exactly $k$ cycles 
of length at least $6$. 
\end{lem}
\noindent
\textbf{Proof of Lemma~\ref{lem:M-2-factor}.}~Choose $k+1$ disjoint $M$-cycles $C_{1}, \dots, C_{k+1}$ in $G$ 
so that 
\begin{align}
\label{maximal cycles}
\sum_{i=1}^{k+1}|C_{i}| \textup{ is as large as possible.} 
\end{align}
We may assume that $V(G) \setminus \bigcup_{i=1}^{k+1}V(C_{i}) \neq \emptyset$.
Let $H = G- \bigcup_{i=1}^{k+1}V(C_{i})$.

By (\ref{maximal cycles}), 
every edge in $M_{H}$ is not insertible in each $M$-cycle $C_{i}$, 
and hence 
by Lemma~\ref{lem:insertion edge}(\ref{insertion edge in cycle}), 
$e_{G}(\{x, y\}, C_{i}) \le \frac{|C_{i}|}{2}$ holds for $xy \in M_{H}$ and $1 \le i \le k + 1$. 
Suppose that $H$ is not connected. 
Let $x_{1}y_{1}, x_{2}y_{2}$ be edges of $M$ in two different components of $H$ 
with orders $2n_{1}, 2n_{2}$ respectively, 
and let $n_{0} = \sum_{i=1}^{k+1}\frac{|C_{i}|}{2}$. 
Then $d_{G}(x_{h}) + d_{G}(y_{h}) \le n_{0} + 2n_{h}$ ($h \in \{1, 2\}$), 
and so 
\begin{align*}
d_{G}(x_{1}) + d_{G}(y_{1}) + d_{G}(x_{2}) + d_{G}(y_{2}) \le 2 (n_{0} + n_{1} + n_{2}) \le 2n < 2n + 4, 
\end{align*}
contrary to hypothesis. 
Thus $H$ is connected.

Consider 
a maximal $M$-path $P_{0}[x, y]$ in $H$, 
and apply Lemma~\ref{lem:crossing} 
with 
$s = k + 1$ and $(D_{1}, \dots, D_{s}, P_{0}) = (C_{1}, \dots, C_{k+1}, P_{0})$. 
If (\ref{P0 cup Di}) holds, then 
replacing $D_{i}$ by $D_{i}'$ contradicts (\ref{maximal cycles}). 
If (\ref{degree sum of x and y in outside}) holds, then 
$e_{G}(\{x , y\}, C_{h}) \ge |C_{h}|/2 + 1$ for some $h$ with $1 \le h \le k+1$, 
that is, $P_{0}$ is insertible in $C_{h}$ 
(by Lemma~\ref{lem:insertion edge}(\ref{insertion edge in cycle})), 
which contradicts (\ref{maximal cycles}) again. 
Thus (\ref{P0}) holds. 
Since $H$ is connected and $P_{0}$ is maximal, 
this implies that 
$H$ has an $M_{H}$-Hamilton cycle if $|H| \ge 4$.

By the degree sum condition, 
$G$ is connected and so 
$E_{G}(H, C_{i}) \neq \emptyset$ for some $i$ with $1 \le i \le k+1$. 
We may assume that $i=1$. 
Since $H$ has an $M_{H}$-Hamilton cycle if $|H| \ge 4$, 
it follows that 
$G[V(H \cup C_{1})]$ has an $M_{H \cup C_{1}}$-Hamilton path $P_{0}'$. 
We now apply Lemma~\ref{lem:crossing} 
with $s = k$ and $(D_{1}, \dots, D_{s}, P_{0}) = (C_{2}, \dots, C_{k+1}, P_{0}')$. 
Note that $|P_{0}'| \ge 6$. 
If 
(\ref{P0}) or (\ref{P0 cup Di}) holds, 
then this contradicts (\ref{maximal cycles}). 
Thus (\ref{degree sum of x and y in outside}) holds. 
This together with Lemma~\ref{lem:insertion edge}(\ref{insertion edge in cycle}) implies 
$P_{0}'$ is insertible in $C_{h}$ for some $h$ with $2 \le h \le k + 1$. 
We name the $M_{P_{0}' \cup C_{h}}$-Hamilton cycle of $G[V(P_{0}' \cup C_{h})]$ as $C_{h}'$, 
and then 
$C_{2} \cup \cdots \cup C_{h-1} \cup C_{h}' \cup C_{h+1} \cup \dots \cup C_{k+1}$ 
forms an $M$-2-factor with exactly $k$ cycles of length at least $6$. 
\qed

Now we are ready to prove Theorem~\ref{thm:from k+1 disjoint cycles to M-2-factor with k cycles}.

\bigskip
\noindent
\textbf{Proof of Theorem~\ref{thm:from k+1 disjoint cycles to M-2-factor with k cycles}.}~By Lemma~\ref{lem:M-2-factor}, 
we may assume that 
$G$ has an $M$-$2$-factor
with exactly $k+1$ cycles $C_{1},\ldots,C_{k+1}$ of length at least $6$,
but
$G$ has no $M$-$2$-factor with exactly $k$ cycles of length at least $6$ 
(since otherwise, the result holds). 
Then, the following fact holds.

\begin{fact}
\label{desired 2-factor}
Let $I \subseteq \{1, 2, \dots, k + 1\}$ with $|I| \ge 2$, and let $C = \bigcup_{i \in I}C_{i}$. 
Then 
$G[V(C)]$ does not have an $M_{C}$-2-factor with exactly $|I| - 1$ cycles of length at least $6$. 
\end{fact}

Let $X$ and $Y$ be partite sets of $G$, 
and we may assume that $xx^{+} \in E(C_{i}) \cap M$ 
for each cycle $C_{i}$ and $x \in V(C_{i}) \cap X$.

\begin{claim}\label{Claim:parallel}
For any two cycles $C_{i}$ and $C_{j}$ with $i \neq j$,
there exist no four vertices $x_{i} \in V(C_{i}) \cap X$, $x_{j} \in V(C_{j}) \cap X$,
$y_{i} \in V(C_{i}) \cap Y$ and $y_{j} \in V(C_{j}) \cap Y$
such that $\{x_{i}y_{j}, y_{i}x_{j}\} \subseteq E(G)$.
\end{claim}
\proof
Suppose that there exist two cycles and four vertices as specified in the claim. 
We may assume that $i = 1$ and $j = 2$. 
Choose $C_{1}, C_{2}, x_{1}, y_{1}, x_{2}$ and $ y_{2}$
so that 
\begin{align}
\label{choice}
\textup{
$x_{1}y_{1} \in E(C_{1}) \setminus M$ or $x_{2}y_{2} \in E(C_{2}) \setminus M$ if possible.
} 
\end{align}

\begin{figure}[h]
\begin{center}
\includegraphics[scale=1.00,clip]{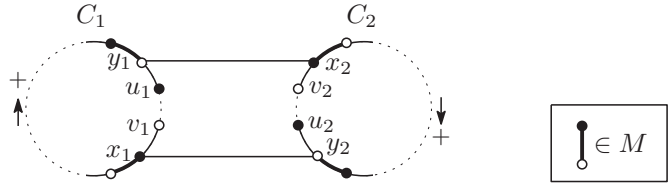}
\caption{The vertices $x_{h}, y_{h}, u_{h}, v_{h}$ ($h \in \{1, 2\}$)}
\label{twoedges}
\end{center}
\end{figure}

Let $u_{1}=y_{1}^{+}$, $v_{1}=x_{1}^{-}$, $u_{2}=y_{2}^{+}$ and $v_{2}=x_{2}^{-}$ 
(see Figure~\ref{twoedges}). 
We consider the $M_{C_{1} \cup C_{2}}$-Hamilton paths 
$P_{0} = u_{1} \ora{C_{1}} y_{1} x_{2} \ora{C_{2}} v_{2}$ 
and 
$P_{0}' = v_{1} \ola{C_{1}} x_{1} y_{2} \ola{C_{2}} u_{2}$ in 
$G[V(C_{1} \cup C_{2})]$ 
and 
we apply Lemma~\ref{lem:crossing} 
with 
$s = k - 1$, 
$(D_{1}, \dots, D_{s}, P_{0}) = (C_{3}, \dots, C_{k+1}, P_{0})$ 
and 
$(D_{1}, \dots, D_{s}, P_{0}) = (C_{3}, \dots, C_{k+1}, P_{0}')$, respectively. 
Then, 
either (\ref{P0}), (\ref{P0 cup Di}) or (\ref{degree sum of x and y in outside}) holds 
for $(C_{3}, \dots, C_{k+1}, P_{0})$ and $(C_{3}, \dots, C_{k+1}, P_{0}')$, respectively.
Combining this with Fact~\ref{desired 2-factor}, 
we see that (\ref{degree sum of x and y in outside}) holds for each case. 
Let $H_{1 2} = G - V(C_{1} \cup C_{2})$. 
Then, we get 
\begin{align*}
e_{G}(\{u_{1}, v_{1}\}, H_{12}) 
+ e_{G}(\{u_{2}, v_{2}\}, H_{12}) 
&
= 
e_{G}(\{u_{1}, v_{2}\}, H_{12})  
+ 
e_{G}(\{v_{1}, u_{2}\}, H_{12}) \\
&
\ge \Big( \frac{|H_{12}|}{2}+1 \Big) + \Big( \frac{|H_{12}|}{2}+1 \Big) = |H_{12}| + 2. 
\end{align*}
This implies that 
$$e_{G}(\{u_{1}, v_{1}\}, H_{12}) \ge \frac{|H_{12}|}{2}+1 
\text{ or }
e_{G}(\{u_{2}, v_{2}\}, H_{12}) \ge \frac{|H_{12}|}{2}+1.$$

\noindent
Without loss of generality,
we may assume that
$e_{G}(\{u_{1}, v_{1}\}, H_{12}) \ge \frac{|H_{12}|}{2}+1$.
Then
there exists a cycle $C_{h}$ with $3 \le h \le k+1$ 
such that $e_{G}(\{u_{1}, v_{1}\}, C_{h}) \ge \frac{|C_{h}|}{2}+1$. 
Hence by Lemma~\ref{lem:insertion edge}(\ref{insertion edge in cycle}), 
there is an insertion edge $b_{h}b_{h}^{+}$ of $u_{1} \ora{C_{1}} v_{1}$ in $C_{h}$ 
such that $u_{1}b_{h}, v_{1}b_{h}^{+} \in E(G)$. 
This also implies that 
$x_{1}y_{1} \in E(C_{1}) \setminus M$ or $x_{2}y_{2} \in E(C_{2}) \setminus M$ 
(since otherwise, 
replacing $(C_{1}, C_{2}, x_{1}, y_{1}, x_{2}, y_{2})$ 
with $(C_{1}, C_{h}, u_{1}, v_{1}, b_{h}^{+}, b_{h})$, 
this contradicts the choice (\ref{choice})).

\begin{figure}[h]
\begin{center}
\includegraphics[scale=1.00,clip]{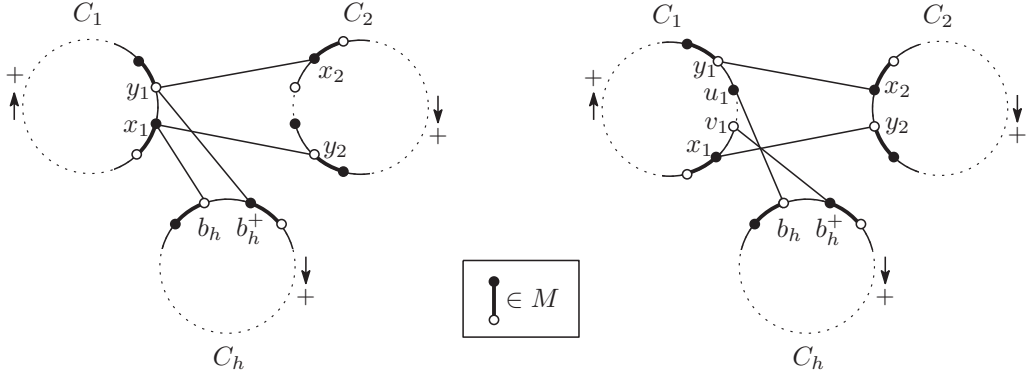}
\caption{Two cases for Claim~\ref{Claim:parallel}}
\label{Claim1}
\end{center}
\end{figure}

If 
$x_{1}y_{1} \in E(C_{1}) \setminus M$, 
that is,
$u_{1}=x_{1}$, 
then 
$x_{1} \ora{C_{1}} y_{1}$ is insertible in $C_{h}$, 
which contradicts Fact~\ref{desired 2-factor} (see the left figure of Figure~\ref{Claim1}).

Thus, $x_{1}y_{1} \notin E(C_{1}) \setminus M$, 
and hence 
$x_{2}y_{2} \in E(C_{2}) \setminus M$,
that is,
$u_{2}=x_{2}$. 
Then, 
$x_{1} \ora{C_{1}} y_{1}x_{2} \ora{C_{2}} y_{2} x_{1}$ 
and 
$u_{1} \ora{C_{1}} v_{1} b_{h}^{+} \ora{C_{h}} b_{h} u_{1}$ 
form an $M_{C_{1} \cup C_{2} \cup C_{h}}$-2-factor with two cycles of length at least $6$ in 
$G[V(C_{1} \cup C_{2} \cup C_{h})]$, 
which contradicts Fact~\ref{desired 2-factor} (see the right figure of Figure~\ref{Claim1}). 
\qed

Since $2n > 6(k+1)$, 
we may assume that $|C_{1}| \ge 8$. 
Note that $E_{G}(C_{1}, G- C_{1}) \neq \emptyset$ 
because $\sigma_{1, 1}(G) \ge n + 2$ implies $G$ is connected and 
$k+1 \ge 2$, 
and hence 
we may also assume
that
there exist two vertices
$y_{1} \in V(C_{1}) \cap Y$ and
$x_{2} \in V(C_{2}) \cap X$ such that $y_{1}x_{2} \in E(G)$. 
Let $x_{1} \in V(C_{1})$ and $y_{2} \in V(C_{2})$ 
with 
$x_{1}y_{1} \in E(C_{1}) \cap M$ and 
$x_{2}y_{2} \in E(C_{2}) \cap M$ 
(i.e., $x_{1} = y_{1}^{-}$ and $y_{2} = x_{2}^{+}$). 
Let $u_{1}=y_{1}^{+}$, $v_{1} = x_{1}^{-}$ and $v_{2}=x_{2}^{-}$.

\begin{figure}[h]
\begin{minipage}{.48\linewidth}
\begin{center}
\includegraphics[scale=1.00,clip]{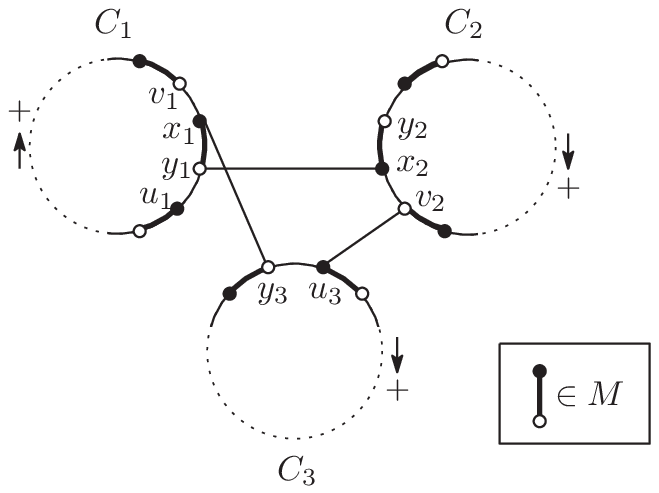}
\caption{$\{x_{1}y_{3}, v_{2}u_{3} \} \subseteq E(G)$}
\label{b3b3+}
\end{center}
\end{minipage}
\begin{minipage}{.48\linewidth}
\begin{center}
\includegraphics[scale=1.00,clip]{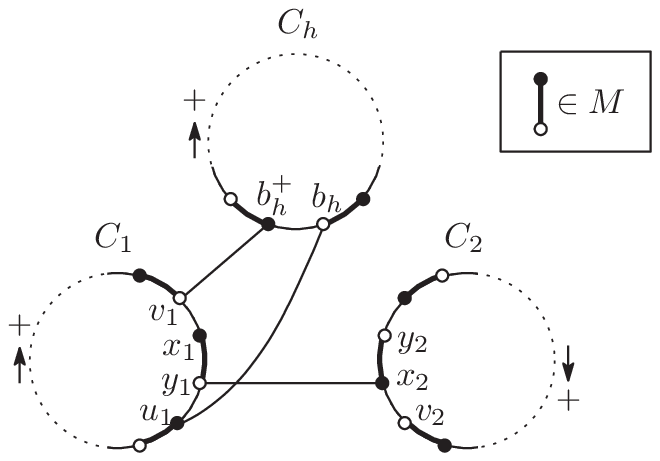}
\caption{$e_{G}(\{u_{1}, v_{1}\}, C_{h}) \ge \frac{|C_{h}|}{2}+1$}
\label{Ch}
\end{center}
\end{minipage}
\end{figure}

By Claim \ref{Claim:parallel},
$E_{G}(x_{1}, C_{2}) = E_{G}(v_{2}, C_{1}) = \emptyset$.
In particular,
$x_{1}v_{2} \not\in E(G)$.
Let $H_{12}=G-V(C_{1} \cup C_{2})$.
Then we obtain
\begin{align*}
e_{G}(\{x_{1}, v_{2}\}, H_{12})
&
= 
\big( d_{G}(x_{1})+d_{G}(v_{2}) \big)
- d_{G[V(C_{1})]}(x_{1}) - e_{G}(x_{1}, C_{2}) 
- e_{G}(v_{2}, C_{1}) - d_{G[V(C_{2})]}(v_{2}) \\
&
\ge 
( n+2 ) 
- |C_{1}|/2 - 0 - 0 - |C_{2}|/2 = |H_{12}|/2+2. 
\end{align*}
Therefore $k \ge 2$ and there exists 
a cycle $C_{h}$ with $3 \le h \le k + 1$ 
such that
$e_{G}(\{x_{1}, v_{2}\}, C_{h}) \ge |C_{h}|/2+1$.
Without loss of generality,
we may assume that $h=3$.
Then by Lemma~\ref{lem:insertion edge}(\ref{insertion edge in cycle}), 
there exists an edge $y_{3}u_{3}$ 
with $u_{3} \in V(C_{3}) \cap X$ and $y_{3} \in V(C_{3}) \cap Y$ 
in $E(C_{3}) \setminus M$ 
such that 
\begin{align}
\label{x1y3, v2u3}
\{x_{1}y_{3}, v_{2}u_{3}\} \subseteq E(G) 
\end{align}
(see Figure~\ref{b3b3+}).

Now consider
the $M$-path $P_{0} = u_{1} \ora{C_{1}} v_{1}$ 
and
the $M$-cycle $D_{0} = x_{1}y_{1}x_{2} \ora{C_{2}} v_{2} u_{3} \ora{C_{3}} y_{3} x_{1}$ 
in $G[V(C_{1} \cup C_{2} \cup C_{3})]$. 
Recall that $|C_{1}| \ge 8$, that is, $|P_{0}| \ge 6$. 
We apply Lemma~\ref{lem:crossing} 
with 
$s = k - 1$ and 
$(D_{1}, \dots, D_{s}, P_{0}) = (D_{0}, C_{4}, \dots, C_{k+1}, P_{0})$.
Then, 
either (\ref{P0}), (\ref{P0 cup Di}) or (\ref{degree sum of x and y in outside}) holds 
for $(D_{0}, C_{4}, \dots, C_{k+1}, P_{0})$. 
Combining this with Fact~\ref{desired 2-factor}, 
we see that (\ref{degree sum of x and y in outside}) holds. 
Let $H_{0}=G-P_{0}$.
Then
we get 
$e_{G}(\{u_{1}, v_{1}\}, H_{0}) 
\ge |H_{0}|/2 +1$. 
Hence either 
(a) there exists a cycle $C_{h}$ with $4 \le h \le k$ 
such that $e_{G}(\{u_{1}, v_{1}\}, C_{h}) \ge \frac{|C_{h}|}{2}+1$, 
or
(b) $e_{G}(\{u_{1}, v_{1}\}, D_{0}) \ge \frac{|D_{0}|}{2}+1$. 
If (a) holds, 
then by Lemma~\ref{lem:insertion edge}(\ref{insertion edge in cycle}), 
there exists an edge $b_{h}b_{h}^{+}$ in $E(C_{h})$ 
such that $u_{1}b_{h}, v_{1}b_{h}^{+} \in E(G)$,
which contradicts Claim \ref{Claim:parallel}
(see Figure~\ref{Ch}). 
Thus 
$e_{G}(\{u_{1}, v_{1}\}, D_{0}) \ge \frac{|D_{0}|}{2}+1$. 
Since 
$x_{1}y_{3} \in E_{G}(C_{1}, C_{3})$ 
and 
$y_{1}x_{2} \in E_{G}(C_{1}, C_{2})$, 
it follows from Claim \ref{Claim:parallel} that 
$E_{G}(v_{1}, C_{3}) = \emptyset$
and
$E_{G}(u_{1}, C_{2}) = \emptyset$.
Therefore, 
the inequality $e_{G}(\{u_{1}, v_{1}\}, D_{0}) \ge \frac{|D_{0}|}{2}+1$ implies that 
\begin{align}
\label{neighbors of u1 and v1}
V(C_{3}) \cap Y \subseteq  N_{G}(u_{1}) 
\textup{ and }
V(C_{2}) \cap X \subseteq  N_{G}(v_{1}).
\end{align}

Let $x_{3} \in V(C_{3})$ with $x_{3}y_{3} \in E(C_{3}) \cap M$.

\begin{figure}[h]
\begin{minipage}{.48\linewidth}
\begin{center}
\includegraphics[scale=1.00,clip]{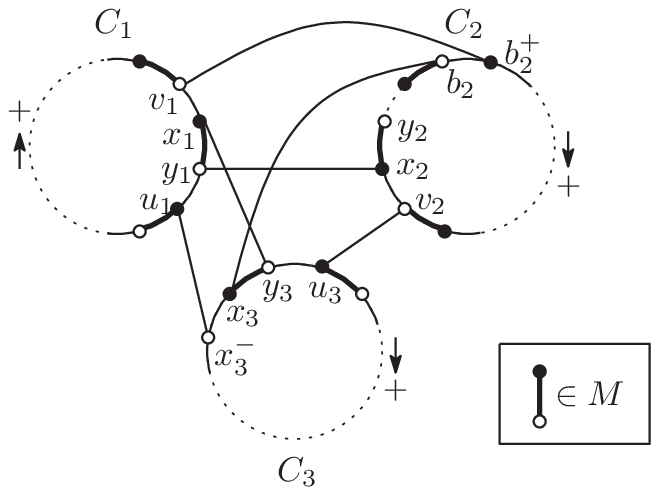}
\caption{Claim \ref{x3}}
\label{figx3}
\end{center}
\end{minipage}
\begin{minipage}{.48\linewidth}
\begin{center}
\includegraphics[scale=1.00,clip]{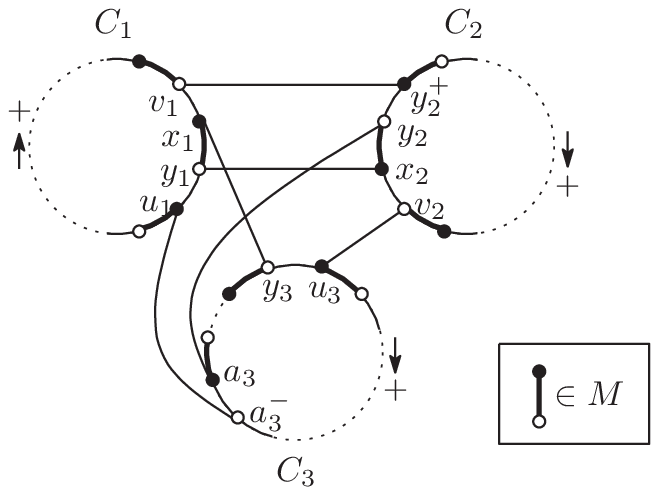}
\caption{Claim \ref{y2}}
\label{v3}
\end{center}\end{minipage}
\end{figure}

\begin{claim}
\label{x3}
$E_{G}(x_{3}, C_{2} - \{v_{2}\}) = \emptyset$.
\end{claim}
\proof
Suppose that 
there exists a vertex $b_{2} \in V(C_{2}) \setminus \{v_{2}\}$
with $x_{3}b_{2} \in E(G)$.
Recall that by (\ref{x1y3, v2u3}), $x_{1}y_{3}, v_{2}u_{3} \in E(G)$. 
Recall also that by (\ref{neighbors of u1 and v1}), $u_{1}x_{3}^{-}, v_{1}b_{2}^{+} \in E(G)$. 
Then 
$u_{1} \ora{C_{1}} v_{1} b_{2}^{+} \ora{C_{2}} v_{2} u_{3} \ora{C_{3}} x_{3}^{-} u_{1}$ 
and 
$x_{1}y_{1}x_{2} \ora{C_{2}} b_{2} x_{3} y_{3} x_{1}$ 
form an $M_{C_{1} \cup C_{2} \cup C_{3}}$-2-factor with two cycles of length at least $6$ in 
$G[V(C_{1} \cup C_{2} \cup C_{3})]$, 
which contradicts Fact~\ref{desired 2-factor} 
(see Figure~\ref{figx3}). 
\qed

\begin{claim}
\label{y2}
$E_{G}(y_{2}, C_{3} - \{u_{3}\}) = \emptyset$.
\end{claim}
\proof
Suppose that 
there exists a vertex $a_{3} \in V(C_{3}) \setminus \{u_{3}\}$
with $y_{2}a_{3} \in E(G)$. 
Recall that by (\ref{x1y3, v2u3}), $x_{1}y_{3}, v_{2}u_{3} \in E(G)$. 
Recall also that by (\ref{neighbors of u1 and v1}), $u_{1}a_{3}^{-}, v_{1}y_{2}^{+} \in E(G)$. 
Then 
$u_{1} \ora{C_{1}} v_{1} y_{2}^{+} \ora{C_{2}} v_{2} u_{3} \ora{C_{3}} a_{3}^{-} u_{1}$ 
and 
$x_{1} y_{1} x_{2} y_{2} a_{3} \ora{C_{3}} y_{3} x_{1}$ 
form an $M_{C_{1} \cup C_{2} \cup C_{3}}$-2-factor with two cycles of length at least $6$ in 
$G[V(C_{1} \cup C_{2} \cup C_{3})]$, 
which contradicts Fact~\ref{desired 2-factor} 
(see Figure~\ref{v3}). 
\qed

Let $H_{23}=G-V(C_{2} \cup C_{3})$. 
Since $v_{2}u_{3} \in E_{G}(C_{2}, C_{3})$, 
it follows from Claim~\ref{Claim:parallel} 
that 
$x_{2}y_{3} \not\in E(G)$ and
\begin{align*}
e_{G}(\{x_{2}, y_{3}\}, H_{23})
&
= 
\big( d_{G}(x_{2})+d_{G}(y_{3}) \big) 
- d_{G[V(C_{2})]}(x_{2}) - e_{G}(x_{2}, C_{3}) 
- e_{G}(y_{3}, C_{2}) - d_{G[V(C_{3})]}(y_{3}) \\
&
\ge ( n + 2 ) - |C_{2}|/2 - 0 - 0 - |C_{3}|/2
= |H_{23}|/2+2.
\end{align*}

\noindent
On the other hand, 
by Claims~\ref{x3} and \ref{y2},
we have
$x_{3}y_{2} \not\in E(G)$ and
\begin{align*}
e_{G}(\{x_{3}, y_{2}\}, H_{23}) 
&
= 
\big( d_{G}(x_{3})+d_{G}(y_{2}) \big) 
-e_{G}(x_{3}, C_{2}) - d_{G[V(C_{3})]}(x_{3}) 
-d_{G[V(C_{2})]}(y_{2}) - e_{G}(y_{2}, C_{3})\\
&
\ge 
(n + 2) - 1 - |C_{3}|/2 - |C_{2}|/2 - 1
= |H_{23}|/2.
\end{align*}

\noindent
Therefore we deduce
$$e_{G}(\{x_{2}, y_{2}, x_{3}, y_{3}\}, H_{23})  \ge |H_{23}|+2.$$
This implies that 
$$
e_{G}(\{x_{2}, y_{2}\}, H_{23}) \ge \frac{|H_{23}|}{2}+1 
\text{ or }
e_{G}(\{x_{3}, y_{3}\}, H_{23}) \ge \frac{|H_{23}|}{2}+1.
$$

\noindent
Without loss of generality,
we may assume that
$e_{G}(\{x_{2}, y_{2}\}, H_{23}) \ge \frac{|H_{23}|}{2}+1$.
Then
there exists a cycle $C_{h}$ with $h = 1$ or $4 \le h \le k+1$ 
such that $e_{G}(\{x_{2}, y_{2}\}, C_{h}) \ge \frac{|C_{h}|}{2}+1$. 
But this clearly leads to a contradiction with Claim \ref{Claim:parallel}, 
and this completes the proof of  Theorem~\ref{thm:from k+1 disjoint cycles to M-2-factor with k cycles}. 
\qed

\section{Related problems}
\label{sec:related problems}

To conclude the paper, we mention related problems. 
In order to prove our main result (Theorem~\ref{thm: sigma1,1 for directed 2-factor with k cycles}), 
we have considered $k$ disjoint $M$-cycles of length at least $6$ in bipartite graphs (Step~1) 
and have given Theorem~\ref{thm:small k disjoint M-cycles} in Section~\ref{sec:proof of main}. 
However, 
in the proof of Theorem~\ref{thm: sigma1,1 for directed 2-factor with k cycles}, 
we do not use ``each $M$-cycle has length at most $8$'' in the conclusion of Theorem~\ref{thm:small k disjoint M-cycles} 
(see the proof in Section~\ref{sec:proof of main}). 
Therefore, 
one may consider the problem whether 
the degree condition can be 
weakened if we drop 
the condition ``length at most $8$'' in Theorem~\ref{thm:small k disjoint M-cycles}. 
Note that a much weaker degree condition than the one of Theorem~\ref{Thm:BCFGL1997} 
guarantees the existence of $k$ disjoint cycles in simple undirected graphs.

\begin{Thm}[Enomoto \cite{Enomoto1998}, Wang \cite{Wang1999-2}]
\label{Thm:Enomoto1998, Wang1999-2}
Let $k$ be a positive integer, 
and let $G$ be a graph of order at least $3k$. 
If $\sigma_{2}(G) \ge 4k - 1$, 
then $G$ contains $k$ disjoint cycles. 
\end{Thm}

Considering the relation between simple undirected graphs, digraphs and bipartite graphs 
(Remarks~\ref{rem:D_{G}} and \ref{rem:from D to G}), 
we can consider the following problem as a digraph version of Theorem~\ref{Thm:Enomoto1998, Wang1999-2}.

\begin{pro}
\label{pro:}
Let $k$ be a positive integer, 
and let $D$ be a digraph of order $n \ge 3k$. 
Is it true that, if $d_{D}^{+}(u) + d_{D}^{-}(v) \ge 4k-1$ 
for every two distinct vertices $u$ and $v$ with $(u, v) \notin A(D)$, 
then $D$ contains $k$ disjoint directed cycles of length at least $3$? 
\end{pro}

As another related problem, 
the following ``minimum out-degree'' condition is conjectured by Bermond and Thomassen (1981) 
for the existence of just $k$ disjoint directed cycles.

\begin{Con}[Bermond, Thomassen \cite{BT1981}]
\label{Con:BT1981}
Let $k$ be a positive integer, 
and let $D$ be a digraph. 
If every vertex has out-degree at least $2k-1$, 
then $D$ contains $k$ disjoint directed cycles. 
\end{Con}

The case $k = 1$ of this conjecture is an easy problem, 
and 
the cases $k = 2$ and $k = 3$ are proved 
in \cite{Thomassen1983} and \cite{LPS2009}, respectively. 
Alon \cite{Alon1996} proved that the conclusion holds 
if every vertex has out-degree at least $64k$, 
but the conjecture as stated remains open.



\end{document}